\documentclass{article}
\usepackage[arrow, matrix, curve]{xy}
\usepackage{amsthm}
\usepackage{hyperref}
\input{amssym.tex}
\newtheorem{thm}{Theorem}
\newtheorem{lem}{Lemma}
\newtheorem{cor} {Corollary}

\newtheorem{obs} {Observation}
\newtheorem{pro} {Proposition}
\newtheorem{df} {Definition}

\newcommand{\ibd}{\partial_\infty\widetilde{M}}

\newcommand{\ot}{\otimes_{{\Bbb Z}G}{\Bbb Z}}

\newcommand{\ib}{\partial_\infty\widetilde{M}}

 \newcommand{\om}{\overline{M}}
	\newenvironment{pf}{{\it Proof:}\quad}{\hfill$QED$}

\begin{document}

	\title{Group homology and ideal fundamental cycles}
	\author{Thilo Kuessner}
	\date{}
	\maketitle
	\begin{abstract} We prove that the (homological version of the) generalized Goncharov invariant of locally symmetric spaces determines their generalized Neumann-Yang invariant.
\noindent
		\end{abstract}

		\section{Introduction}
Let $M$ be an odd-dimensional 
locally symmetric space of noncompact type, with $vol\left(M\right)<\infty$, which is either compact or of rank one. 
We compare two invariants: the generalized Goncharov invariant and the generalized Neumann-Yang invariant, and we show that the former determines the latter. 

In \cite{ny}, Walter D.\ Neumann and Jun Yang constructed an invariant of finite-volume hyperbolic 3-manifolds which lives in the Bloch group ${\mathcal{B}}\left({\bf C}\right)$, and from which volume and Chern-Simons invariant can be recaptured. On the 
other hand, Alexander Goncharov \cite{gon} constructed an invariant of odd-dimensional hyperbolic manifolds and representations of their fundamental groups. (In particular he considered the half-spinor representations.)
This invariant 
lives in $K_*\left(\overline{\bf Q}\right)\otimes {\bf Q}$ and the volume can be recaptured by means of the Borel regulator. 
In \cite{ku}, we generalized Goncharov's construction to finite-volume locally symmetric spaces of noncompact type $\Gamma\backslash G/K$
(either closed or locally rank one symmetric) and representations $\rho:G\rightarrow GL\left(N,{\bf C}\right)$.
We constructed an invariant $\overline{\gamma}\left(M\right)\in H_d\left(GL\left(N,\overline{\bf Q}\right);{\bf Q}\right)$ and a projection $\pi: H_d\left(GL\left(\overline{\bf Q}\right);{\bf Q}\right)\rightarrow K_d\left(\overline{\bf Q}\right)\otimes{\bf Q}$ yielding an invariant $\gamma\left(M\right):=\pi\left(\overline{\gamma}\left(M\right)\right)\in K_d\left(\overline{\bf Q}\right)\otimes{\bf Q}$,
and we analyzed in which cases the
obtained invariant $\gamma\left(M\right)$ is nontrivial resp.\ trivial. (The nontriviality of $\gamma\left(M\right)$ turned out to depend only on the representation $\rho:G\rightarrow GL\left(N,{\bf C}\right)$. Again, if $\gamma\left(M\right)\not=0$, then $vol\left(M\right)$ can be recaptured from $\gamma\left(M\right)$ by application of the Borel regulator.)

Here we show that the natural evaluation map $ev$ from $H_*\left(GL\left(N,{\bf C}\right)\right)$ to the generalized pre-Bloch groups maps $\overline{\gamma}\left(M\right)$ to the natural generalization of the Neumann-Yang invariant $\beta\left(M\right)$. The latter invariant is defined using ideal fundamental cycles. 

{\bf Theorem 1.} {\em Let $\overline{M}$ be a compact, orientable, connected manifold with boundary such that $M:=\overline{M}-\partial \overline{M}=\Gamma\backslash G/K$ is a finite-volume, locally rank one symmetric space of noncompact type.
Let $\rho:G\rightarrow GL\left(N,{\Bbb C}\right)$ be a representation.
Then} $$ev_*\left(\overline{\gamma}\left(\om\right)\right)=\beta_\rho\left(\om\right).$$ For
$N=2$, 
we can derive the following consequence:

{\bf Corollary 2.} {\em If $M^3$ is hyperbolic of finite volume, then Suslin's homomorphism
$$K_3\left({\bf C}\right)\otimes{\bf Q}\rightarrow {\mathcal{B}}\left({\bf C}\right)\otimes{\bf Q}$$ maps the Goncharov invariant $\gamma\left(M\right)$ to the Neumann-Yang invariant $\beta\left(M\right)\otimes 1$.}\\
\\
Although Corollary 2 seems to be well-known to the experts, at least for the case of closed manifolds, we could not locate a reference. 

This paper makes essential use of the definition of the Goncharov invariant for cusped manifolds (\cite[Section 4]{ku} and as such 
is a continuation and extension of \cite{ku}. We decided to publish it separately because \cite{ku} already has 
a considerable length. The paper is organized as follows: Section 2 states the definitions of the Goncharov invariant 
and the (generalized) Neumann-Yang invariant. 
Section 3 provides a chain map $\hat{\Psi}$ from relative chains to ideal chains which, in particular, sends relative fundamental cycles to ideal fundamental cycles.
This is used in Section 4 to prove Theorem 1.
\section{Definitions}
\subsection{Goncharov invariant}
{\bf Group homology.} For a group G, its classifying space (with respect to the discrete topology on $G$) is the simplicial set $BG$ defined as follows:

- the set of k-simplices
of $BG$ is $$S_k^{simp}\left(BG\right)=\left\{\left(g_1,\ldots,g_k\right): g_1,\ldots,g_k\in G\right\},$$

- the boundary operators are defined by

$\partial_0\left(g_1,\ldots,g_k\right)=\left(g_2,\ldots,g_k\right),$

$ \partial_i\left(g_1,\ldots,g_k\right)=\left(g_1,\ldots,g_ig_{i+1},\ldots,g_k\right)$ for $i=1,\ldots,k-1$,

$\partial_k\left(g_1,\ldots,g_k\right)=\left(g_1,\ldots,g_{k-1}\right)$.

- the degeneracy maps are defined by $s_j\left(g_1,\ldots,g_k\right)=\left(g_1,\ldots,g_j,1,g_{j+1},\ldots,g_k\right)$.\\
$C_k\left(BG\right)$ is the free abelian group freely generated by $S_k\left(BG\right)$. Let $\partial:C_k\left(BG\right)\rightarrow C_{k-1}\left(BG\right)$ be the linear extension 
of $\partial=\sum_{i=0}^k\left(-1\right)^i\partial_i$.
For an abelian group $R$, the group homology $H_*\left(G;R\right)=H_*^{simp}\left(BG;R\right)$ is the homology of $\left(C_*\left(BG\right)\otimes_{\bf Z}R,\partial\otimes 1\right)$.

If $\Gamma=\pi_1M$ for an aspherical space $M$, then $H_*^{simp}\left(B\Gamma;R\right)=H_*\left(M;R\right)$. In particular, this is the case if $M$ is a locally symmetric space of noncompact type.\\
\\
{\bf Homological Goncharov invariant of closed manifolds.} Let $M^d=\Gamma\backslash G/K$ be a closed, orientable, connected, locally symmetric space of noncompact type, and $\rho:G\rightarrow SL\left(N,{\bf C}\right)$ 
a representation\footnotemark\footnotetext[1]{Since $G$ is semisimple, any representation $\rho:G\rightarrow GL\left(N;{\bf C}\right)$ has image in $SL\left(N,{\bf C}\right)$.}
. Let $\left[M\right]$ be a generator of $
H_d\left(M;{\bf Z}\right)\simeq {\bf Z}$.
Then the homological
Goncharov invariant of $\left(M,\rho\right)$ is defined as $$\overline{\gamma}\left(M\right):=B\left(\rho j\right)_dEM_d^{-1}\left[M\right]\in H_d^{simp}\left(BSL\left(N,{\bf C}\right);{\bf Z}\right),$$
with $j:\Gamma\rightarrow G$ the inclusion, $B\left(\rho j\right):B\Gamma\rightarrow BSL\left(N,{\bf C}\right)$ induced by $\rho j$ and $EM_*^{-1}:H_*^{simp}\left(B\Gamma;{\bf Z}\right)\rightarrow H_*\left(M;{\bf Z}\right)$ the Eilenberg-MacLane isomorphism (cf.\ \cite{ku}, Section 2.1).

(In \cite{ku}, $\overline{\gamma}\left(M\right)$ was considered as an element of $H_d\left(BGL\left({\bf C}\right);{\bf Q}\right)$. Moreover, in Section 
2.5.\ of \cite{ku} we constructed a projection $H_d\left(BGL\left({\bf C}\right);{\bf Q}\right)\rightarrow PH_d\left(BGL\left({\bf C}\right);{\bf Q}\right)\simeq K_d\left({\bf C}\right)\otimes{\bf Q}$ and defined the {\bf K-theoretic Goncharov invariant} $\gamma\left(M\right)\in K_d\left({\bf C}\right)\otimes{\bf Q}$ as the image of $\overline{\gamma}\left(M\right)$ under this projection.

In this paper we will,
except for \hyperref[suslin]{Corollary \ref*{suslin}}, only 
use the homological invariant $\overline{\gamma}\left(M\right)$.
The same remark applies to cusped manifolds below.)\\
\\
{\bf Homological Goncharov invariant of cusped manifolds.}
We recall some notation from \cite[Section 4.2]{ku}:
For a continuous mapping $:A_1\dot{\cup}\ldots\dot{\cup}A_s\rightarrow X$ we define  the {\bf disjoint cone} $$DCone\left(\cup_{i=1}^sA_i\rightarrow X\right)$$
to be the pushout of the diagram
$$ \begin{xy}
\xymatrix{
A_1\dot{\cup}\ldots\dot{\cup}A_s\ar[r]^i \ar[d]&
X \ar[d]\\
Cone\left(A_1\right)\dot{\cup}\ldots\dot{\cup}Cone\left(A_s\right)\ar[r]&DCone\left(\cup_{i=1}^sA_i\rightarrow X\right)}
\end{xy}$$
If $X$ is a CW-complex and $A_1,\ldots,A_s$ are disjoint sub-CW-complexes, then clearly
$$H_*\left(DCone\left(\cup_{i=1}^sA_i\rightarrow X\right)\right)\cong
H_*\left(Cone\left(\cup_{i=1}^sA_i\rightarrow X\right)\right)=H_*\left(X,\cup_{i=1}^s A_i\right)$$
in degrees $*\ge 2$. Similarly, for a simplicial mapping of simplicial sets 
$B_1\dot{\cup}\ldots\dot{\cup}B_s\rightarrow Y$ we defined $DCone\left(\cup_{i=1}^sB_i\rightarrow Y\right)$ to be the quasi-simplicial set whose $q$-simplices are either $q$-simplices in $Y$ or cones over $q-1$-simplices in some $B_i$, with the natural boundary operator.

\begin{pro}\label{preimage}
 Let
$\overline{M}$ be a compact, orientable, connected $d$-manifold with boundary components
$\partial_1\overline{M},\ldots,\partial_s\overline{M}$ such that
 $M:=\overline{M}-\partial \overline{M}=\Gamma\backslash
         G/K$ is a finite-volume locally rank one
        symmetric space of noncompact type.
Let $\rho:G\rightarrow SL\left(N,{\bf C}\right)$
be a representation
such that $\Gamma_i^\prime:=\rho\left(\Gamma_i\right)$ is unipotent for
$i=1,\ldots,s$, where $\Gamma_i\subset \Gamma$ is the subgroup of $\pi_1\left(M,x\right)$ corresponding to\footnotemark\footnotetext[2]{Our assumptions imply that $\pi_1\partial_iM$ injects into $\pi_1M$, cf.\ \cite{ku}, Observation 1. In particular, if we fix $x_0\in M$ and $x_i\in\partial_iM$ for $i=1,\ldots,s$, then we obtain (using
some path from $x_0$ to $x_i$) {\bf isomorphisms of $\pi_1\left(\partial_iM,x_i\right)$ with subgroups $\Gamma_i$ of $\Gamma$}$=\pi_1\left(M,x_0\right)$, for $i=1,\ldots,s$, and we will assume these isomorphisms to be fixed.

Moreover, since $M$ and all $\partial_iM$ are aspherical (this follows from \cite{ebe}), there 
is an isomorphism $H_*\left(DCone\left(\cup_{i=1}^s B\Gamma_i\rightarrow B\Gamma\right)\right)\rightarrow
H_*\left(DCone\left(\cup_{i=1}^s\partial_i M
\rightarrow M\right)\right)
$, see \cite[Lemma 8]{ku}.} $\pi_1\left(\partial_i M,x_i\right)$.
        Denote $$\left[\overline{M},\partial \overline{M}\right]\in H_d\left(DCone\left(\cup_{i=1}^s\partial_i \overline{M}\rightarrow \overline{M}\right);{\Bbb Z}\right)$$
        the fundamental class of $\overline{M}$. Then $$  B\left(\rho j\right)_dEM_d^{-1}\left[\overline{M},\partial
\overline{M}\right]\in H_d\left(DCone\left(\cup_{i=1}^sB\Gamma_i^\prime\rightarrow BSL\left(N,{\bf C}\right)
\right);{\Bbb Z}\right)$$ has a preimage
$$\overline{\gamma}\left(\overline{M}\right)\in H_d\left(BSL\left(N,{\bf C}\right);{\Bbb Z}\right).$$
\end{pro}
This is proved in \cite[Proposition 2]{ku}.

\subsection{Neumann-Yang invariant}

The following definitions are given in \cite[Section 8]{ny} for hyperbolic manifolds (and
actually in more generality for representations with certain properties), but they generalize in an obvious way to 
manifolds of nonpositive sectional curvature, in particular to 
locally symmetric spaces of noncompact type. (See \cite[Definition 8.1 and Definition 8.6]{bh} for the definition of
the boundary at infinity $\partial_\infty\widetilde{M}$ and of the topology on $\widetilde{M}\cup\partial_\infty\widetilde{M}$.)

\begin{df}\label{prebloch} a) Let $\widetilde{M}$ be a simply connected Riemannian manifold of nonpositive sectional curvature
and $\ib$ its ideal boundary. Let $C_q\left(\ib\right)$ be the free
abelian group 
freely generated by the (q+1)-tuples of points of $\ib$ modulo the relations \\
i) $ \left(z_0,\ldots,z_q\right)=sign\left(\tau\right)\left(z_{\tau\left(0\right)},\ldots,z_{\tau\left(q\right)}\right)$
for any permutation $\tau$ of $\left\{0,\ldots,q\right\}$ and \\
ii) $\left(z_0,\ldots,z_q\right)=0$
whenever $z_i=z_j$ for some $i\not=j$. \\

The operator $\partial: C_q\left(\ib\right)\rightarrow C_{q-1}\left(\ib\right)$ is defined by
$\partial\left(z_0,\ldots,z_q\right)=\sum_{i=0}^q\left(-1\right)^i\left(z_0,\ldots,\hat{z_i},\ldots,z_q\right)$ and linear extension.
Let $G=Isom\left(\widetilde{M}\right)$. We define the generalized pre-Bloch group of $\widetilde{M}$ as
$${\mathcal{P}}_n\left(\widetilde{M}\right):=H_n\left(C_*\left(\ib\right)\ot, \partial \otimes_{{\Bbb Z}G}  id\right).$$
b) 
We define the generalized pre-Bloch groups of ${\Bbb C}$ as $${\mathcal{P}}_n^N\left({\bf C}\right):=
{\mathcal{P}}_n\left(SL\left(N,{\bf C}\right)/SU\left(N\right)\right).$$\end{df}

{\bf Relation to classical Bloch group.} In particular ${\mathcal{P}}_3^2\left({\bf C}\right)={\mathcal{P}}_3\left({\bf H}^3\right)$ is the classical pre-Bloch group ${\mathcal{P}}\left({\bf C}\right)$
(\cite{ny}\cite{su}).

One should note that Neumann and Yang \cite{ny} considered $C_*\left(\partial_\infty {\bf H}^3\right)
$ as the free abelian group generated by tuples of {\em distinct} 
points in $\partial_\infty {\bf H}^3$. We will denote this group by $C_*^{nd}\left(\partial_\infty{\bf H}^3\right)$, the superscript standing for 'nondegenerate'. It is easy to see that the inclusion $C_*^{nd}\left(\partial_\infty{\bf H}^3\right)\ot\rightarrow C_*\left(\partial_\infty{\bf H}^3\right)\ot$ induces an isomorphism in homology. (Here
we consider the action of $G:=PGL\left(2,{\bf C}\right)$
 on $P^1{\bf C}=\partial_\infty{\bf H}^3$.)

The definition of the pre-Bloch group (\cite[Definition 2.1]{ny})
is easily seen to be equivalent to 
${\mathcal{P}}\left({\bf C}\right)=C_3^{nd}\left(\partial_\infty{\bf H}^3\right)\ot/\partial C_4^{nd}\left(\partial_\infty{\bf H}^3\right)\ot$.
This agrees with $H_3\left(C_*^{nd}\left(\partial_\infty{\bf H}^3\right)\ot, \partial \otimes_{{\Bbb Z}G}  id\right)$ because $\partial\otimes 1:C_3^{nd}\left(\partial_\infty{\bf H}^3\right)\ot\rightarrow C_2^{nd}\left(\partial_\infty{\bf H}^3\right)\ot$ is trivial, see \cite[Remark 4.4]{ny}.
\begin{df}\label{crossratio} Let $M$ be a manifold of nonpositive sectional curvature, $\widetilde{M}$ its universal covering with deck group $\Gamma$. Let $\sigma\in \overline{S}_q^{str}\left(M\right)$ 
be a
proper ideal q-simplex (\hyperref[straight]{Definition \ref*{straight}}), $\tilde{\sigma}:\Delta^q\rightarrow\widetilde{M}
\cup\partial_\infty\widetilde{M}$ a lift of $\sigma$ with vertices $\tilde{\sigma}\left(v_0\right),\ldots,\tilde{\sigma}\left(v_q\right)\in \partial_\infty \widetilde{M}$. 
We define the cross-ratio of $\sigma$ by
$$cr\left(\sigma\right) := \left(\tilde{\sigma}\left(v_0\right),\ldots,\tilde{\sigma}\left(v_q\right)\right)\otimes 1\in C_q\left(\ib\right)\otimes_{{\bf Z}\Gamma}{\bf Z}.$$\end{df}
This is well-defined because all lifts of $\sigma$ are in the $\Gamma$-orbit of $\tilde{\sigma}$.

If $G/K$ is a rank one symmetric space of
noncompact type
and $\rho:G\rightarrow SL\left(N,{\bf C}\right)$ a representation, then 
upon conjugation we can assume 
that $\rho$ maps $K$ to $SU\left(N\right)$, inducing a smooth map 
$\rho:G/K\rightarrow SL\left(N,{\bf C}\right)/SU\left(N\right)$ 
and a continuous $\rho$-equivariant map $\rho_\infty:\partial_\infty G/K\rightarrow\partial_\infty SL\left(N,{\bf C}\right)/SU\left(N\right)$. 

In \cite{ny}, the invariant $\beta\left(M\right)$ of hyperbolic 3-manifolds $M$ was defined by means of degree one ideal triangulations. To circumvent the question whether $M$ admits such a triangulation we give an a priori weaker definition and we will show in Section 3 that this one 
agrees with the Neumann-Yang definition. 

The notion of 'proper ideal fundamental cycle' is defined in \hyperref[ifc]{Definition \ref*{ifc}}.


\begin{df}\label{nyinvariante} 
i) Let $\overline{M}$ be a compact, orientable, connected $d$-manifold with boundary such that $M=\overline{M}- \partial \overline{M}$ is a Riemannian manifold of nonpositive sectional curvature. Let $G$ be the isometry group of the universal covering $\widetilde{M}$.

Let $\sum_{i=1}^ra_i\tau_i$ be a proper ideal fundamental cycle of $M$. The Neumann-Yang invariant of $\overline{M}$ is defined by $$\beta\left(\overline{M}\right):=\sum_{i=1}^r a_i \left[cr\left(\tau_i\right)\right]\in 
H_d\left(C_*\left(\ib\right)\ot, \partial \otimes_{{\Bbb Z}G}  id\right)={\mathcal{P}}_d\left(\widetilde{M}\right).$$
ii) If $\widetilde{M}=G/K$ is a rank one symmetric space of noncompact type, and if
$\rho:G\rightarrow SL\left(N,{\bf C}\right)$ is a representation, then we define
$$\beta_\rho\left(\overline{M}\right):=\rho_*\left(\beta\left(\overline{M}\right)\right)\in{\mathcal{P}}_d^N\left({\bf C}\right),$$
where $\rho_*:H_*\left(C_*\left(\partial_\infty G/K\right)\ot\right)\rightarrow H_n\left(C_*\left(\partial_\infty SL\left(N,{\bf C}\right)/SU\left(N\right)\right)\ot\right)$
is induced by the equivariant chain map $\rho_*\left(c_0,\ldots,c_q\right)=\left(\rho_\infty\left(c_0\right),\ldots,\rho_\infty\left(c_q\right)\right)$. \end{df}

We will show in Section 3.4 that $\beta\left(M\right)$ is well defined, that is, independent of the chosen proper ideal fundamental cycle.\\

{\bf Some results of Neumann-Yang.}
Let $\overline{M}$ be a compact, orientable, connected $d$-manifold with boundary such that $M:=\overline{M}-\partial\overline{M}$ is a Riemannian manifold of negative sectional curvature {\em and finite volume}.
In \cite{ny} it was shown that there is an isomorphism $H_d\left(\overline{M},\partial \overline{M}\right)\simeq H_d\left(C_*\left(\partial_\infty\widetilde{M}\right)_\Gamma\right)$ for $\Gamma=\pi_1M$ acting by deck transformations. 
Although this was stated in \cite{ny} for hyperbolic 3-manifolds, the proof in \cite{ny} only 
requires that $M$ is a negatively curved manifold of dimension $d\ge 3$, as we shall indicate now. (This isomorphism will be used in the proof of well-definedness of $\beta\left(M\right)$.)

The reason the argument from \cite{ny} works in this general context is the following elementary fact about discrete isometry groups in negative curvature:\\
{\em If $\widetilde{M}$ is a simply connected manifold of negative sectional
curvature, if
$\Gamma\subset Isom\left(\widetilde{M}\right)$ is discrete of finite covolume, and $x\in\ib$ is not a cusp of $\Gamma$, then $Stab\left(x\right)=\left\{\gamma\in\Gamma:\gamma x=x\right\}$ is an infinite cyclic group.}\\
Proof: If $Stab\left(x\right)$ contained loxodromic isometries $\gamma_1,\gamma_2$ 
with $Fix\left(\gamma_1\right)=\left\{x,y_1\right\}$ and $Fix\left(
\gamma_2\right)=\left\{x,y_2\right\}$ for $y_1\not=y_2$, then $\left[
\gamma_1,\gamma_2\right]\in\Gamma$ would be a parabolic isometry with fixed point $x$, contradicting the assumption that $x$ is not a cusp
of $\Gamma$. Hence all elements of $Stab\left(x\right)$ have a common fixed point $y\not=x$ and fix the geodesic from $x$ to $y$. Thus $Stab\left(x\right)$ is a discrete subgroup of ${\bf R}$, hence infinite cyclic.\\
\\
The isomorphism ${\bf Z}=H_d\left(\overline{M},\partial \overline{M}\right)
\rightarrow H_d\left(C_*\left(\partial_\infty\widetilde{M}\right)_\Gamma\right)$ will 
be defined\footnotemark\footnotetext[3]{${\mathcal{C}}$ is the set of cusps of $\Gamma$. If ${\mathcal{C}}$, then $H_*\left(\Gamma,{\mathcal{C}}\right)$ is understood to be $H_*\left(\Gamma,{\bf Z}\right)$, that is in this case we claim an isomorphism $H_*\left(M;{\bf Z}\right)\simeq H_*\left(\Gamma;{\bf Z}\right)\simeq H_*\left(\Gamma,\partial_\infty\widetilde{M}\right)$.

Recall that, for a G-set $\Omega$, $H_*\left(G,\Omega\right)$ is defined as the homology of $P_*\ot$, where $$\ldots
P_3\rightarrow P_2\rightarrow P_1\rightarrow J\Omega$$
is a ${\Bbb Z}G$-projective resolution of $J\Omega:=ker\left({\Bbb Z}\Omega\rightarrow {\Bbb Z}\right)$, the kernel of the augmentation map.
} as the composition
$$H_d\left(\overline{M},\partial \overline{M}\right)\simeq H_d\left(\Gamma,{\mathcal{C}}\right)\rightarrow H_d\left(\Gamma,\ib\right)
\rightarrow H_d\left(C_*\left(\ib\right)_\Gamma\right).$$

The first isomorphism can be derived by the same argument as in the proof of \cite[Theorem V.1.3]{dd}.
Namely, $C_*\left(\widetilde{\overline{M}}\right)/C_*\left(\partial\widetilde{\overline{M}}\right)$ is acyclic except in degree 1, where 
its cokernel is $K_*:={\Bbb Z}{\mathcal{C}}$. Thus $H_*\left(\Gamma,{\mathcal{C}}\right)
= H_*\left(C_*\left(\widetilde{\overline{M}}\right)/C_*\left(\partial\widetilde{\overline{M}}\right)\ot\right)$.
On the other hand, by definition of relative homology, $H_*\left(C_*\left(\widetilde{\overline{M}}\right)/C_*\left(\partial\widetilde{\overline{M}}\right)\ot\right)=H_*\left(\overline{M},\partial \overline{M}\right)$.

\begin{lem}\label{isomorphism1} $H_d\left(\Gamma,{\mathcal{C}}\right)\rightarrow H_d\left(\Gamma,\ib\right)$
is an isomorphism for $d\ge 3$. 

In particular, if $\mathcal{C}=\emptyset$, then  $H_d\left(\Gamma;{\bf Z}\right)\rightarrow H_d\left(\Gamma,\ib\right)$
is an isomorphism for $d\ge 3$.\end{lem}
\begin{pf} (\cite[Section 3]{ny}) By Shapiro's Lemma, $H_i\left(\Gamma,{\Bbb Z}\left(\ib - {\mathcal{C}}\right)\right)$ is 
the direct sum of $H_i$ of the isotropy groups for the orbits of $\Gamma$ on $\ib - {\mathcal{C}}$.  (This is also true if $\mathcal{C}=\emptyset$.)

The homology of a cyclic group vanishes for $i\ge 2$, hence $H_i\left(\Gamma,{\Bbb Z}\left(\ib - {\mathcal{C}}\right)\right)
=0$ for $i\ge 2$. 
The long exact homology sequence of $J{\mathcal{C}}\rightarrow J\ib\rightarrow {\Bbb Z}
\left(\ib - {\mathcal{C}}\right)$ implies then the desired isomorphism for $d\ge 3$.\end{pf}

\begin{lem}\label{isomorphism2} $H_d\left(\Gamma,\ib\right)\rightarrow H_d\left(C_*\left(\ib\right)_\Gamma\right)$ is an isomorphism.\end{lem}
\begin{pf} (\cite[Proposition 3.2]{ny})
Consider the spectral sequence with $E^1$-term $H_i\left(\Gamma,C_{j-1}\left(\ib\right)\right)$ that
converges to $H_*\left(\Gamma,\ib\right)$. The isotropy group of a pair of points is infinite cyclic or trivial, that is of homological dimension at most 1, hence $E_{p,q}^1=0$ if $q\ge 2$. Thus the only nontrivial $d_1$ is $d_1:H_1\left(\Gamma,C_1\left(\ib
\right)\right)\rightarrow H_1\left(\Gamma,C_0\left(\ib\right)\right)$. In \cite[Lemma 3.3]{ny} it is proven that 
$d_1$ is injective for hyperbolic 3-manifolds. The proof uses only the fact that stabilizers of points and stabilizers of pairs 
of points in $\partial_\infty\widetilde{M}$ are infinitely cyclic or trivial. By the observation above this is true 
for any negatively curved manifold, hence the proof of \cite{ny} applies verbatim. \end{pf}
\begin{cor}\label{discrete} If
$\widetilde{M}$ is a simply connected manifold of negative sectional curvature and dimension $d\ge 3$, and if $\Gamma\subset Isom^+\left(\widetilde{M}\right)$ is a discrete subgroup of finite covolume, then $H_d\left(C_*\left(\ib\right)_\Gamma\right)\simeq{\bf Z}$.\end{cor}

\section{Ideal fundamental cycles}
\subsection{Definition}
We refer to \cite{bh} for basic notions about simply connected manifolds $\widetilde{M}$ of nonpositive sectional curvature, especially for the notion of 
ideal boundary $\partial_\infty\widetilde{M}$ (\cite[Definition 8.1]{bh}) 
and the topology on 
$\widetilde{M}\cup\ibd$ (\cite[Definition 8.5]{bh}). 
For each ordered tuple $\left(v_0,\ldots,v_q\right)\in \left(\widetilde{M}\cup\ibd\right)^{q+1}$, there is a unique straight ideal simplex $str\left(v_0,\ldots,v_q\right)$ with these vertices. We call this simplex a {\em genuine straight simplex} resp.\ a {\em proper ideal straight simplex} if all vertices belong to $\widetilde{M}$ resp.\ all vertices belong to $\partial_\infty\widetilde{M}$. In the sequel, $\Delta^n$ is the standard $n$-simplex and $\Delta^n_0$ its set of vertices.
%
%
\begin{df}\label{straight} 
Let $M$ be a
Riemannian manifold of nonpositive sectional curvature, $\widetilde{M}$ its universal cover, $\Gamma$ a discrete group of
isometries of $\widetilde{M}$ such that $M=\Gamma\backslash\widetilde{M}$, $\pi:\widetilde{M}\rightarrow M$ the canonical projection.\\
a) Let $\overline{C}_n^{str}\left(M \right)$ be the free abelian group generated by 
$$\overline{S}_n^{str}\left(M \right):=\left\{\sigma:\Delta^n-
\Delta^n_0\rightarrow M:\begin{array}{c}\mbox{\ there\ is\ an\ ideal\ straight\ simplex\ }
\tilde{\sigma}:\Delta^n\rightarrow\widetilde{M}\cup\partial_\infty\widetilde{M}\\
\mbox{\ with\ }\pi
\circ\tilde{\sigma}\mid_{\Delta^n-\Delta^n_0}=\sigma\mid_{\Delta^n-\Delta^n_0}\end{array}\right\}.$$
A simplex $\sigma\in \overline{S}_n^{str}\left(M \right)$ 
will be called genuine resp.\ proper ideal if (some, hence each) 
$\widetilde{\sigma}$ is a genuine resp.\ a proper ideal straight simplex in $\widetilde{M}$.\\
b) For $x_0\in M$ define $C_*^{str,x_0}\left(M\right)\subset \overline{C}_*^{str}\left(M \right)$ to be the subcomplex generated by genuine straight simplices with all vertices in $x_0$.\\
c) For $c_0\in\partial_\infty\widetilde{M}$ define $\overline{C}_*^{str,{c_0}}\left(\widetilde{M} \right)\subset \overline{C}_*^{str}\left(M \right)$
to be the subcomplex generated by proper ideal straight simplices such that (some, hence each) $\tilde{\sigma}$ has all vertices in $\Gamma c_0$. 
\end{df}

There is a chain isomorphism $EM:C_*^{simp}\left(B\Gamma\right)\rightarrow C_*^{str,x_0}\left(M\right)$, essentially due to Eilenberg-MacLane (cf.\ Section 2.1 in \cite{ku}).

\begin{df} Let $M$ be an orientable smooth $d$-manifold.
For $\sigma\in \overline{S}_d^{str}\left(M \right)$ and $y\in M$ we 
define $deg_y\left(\sigma\right)=0$ if 
$y\not\in \sigma\left(\Delta^d-\partial\Delta^d\right)$, else we define
$$deg_y\left(\sigma \right):=\sum_{\sigma\left(x\right)=y}
sign\left(det\left(D_x\sigma\right)\right).$$
For a chain $z=\sum_{i=1}^ra_i\sigma_i\in
\overline{C}_d^{str}\left(M \right)$ we define $$deg_y\left(z\right):=\sum_{i=1}^r a_i deg_y\left(\sigma_i\right).$$\end{df}

{\bf Homology invariance of {\em deg} for closed manifolds.} If  $z=\sum_{i=1}^ra_i\sigma_i\in
\overline{C}_d^{str}\left(M \right)$, then it is a consequence of the Sard Lemma that almost 
every $y\in M$ is not contained in $\cup_{i=1}^r \sigma_i\left(\partial\Delta^d\right)$. If $M$ is a closed, orientable, connected $d$-manifold, $\partial z=0$ and all $\sigma_i$ are genuine straight simplices, then for all $y\not\in 
\cup_{i=1}^r \sigma_i\left(\partial \Delta^d\right)$ the 
isomorphism $H_d\left(M,M-\left\{y\right\}\right)
\rightarrow {\bf Z}$ sends the relative homology class of 
$z$ to
$deg_y\left(z\right)$.
In particular 
$deg_y\left(z\right)$ depends 
only on the homology class $\left[z\right]\in H_d\left(M\right)$.

Since $deg_y\left(z\right)$ is the same for almost all $y\in M$ we will occasionally just write $deg\left(z\right)$.

\begin{obs} a) Let $M$ be a closed, orientable, connected Riemannian $d$-manifold of nonpositive sectional
curvature.
Let $z=\sum_{i=1}^r a_i\sigma_i\in \overline{C}_d^{str}\left(M \right)$ 
be a singular cycle consisting of genuine straight simplices. If $z$ represents the fundamental class $\left[M\right]$, 
then $deg_y\left(z\right)=1$ for all $y\not\in \cup_{i=1}^r\sigma_i\left(\partial\Delta^d\right)$.\\
b) 
Let $z=\sum_{i=1}^r \sigma_i\in \overline{C}_d^{str}\left(M \right)$ be an ideal degree one triangulation of a finite volume hyperbolic 3-manifold as in \cite{ny}. Then $deg_y\left(z\right)=1$ for all $y\not\in \cup_{i=1}^r\sigma_i\left(\partial\Delta^d\right)$.\end{obs}

\begin{pf} Smooth closed manifolds are triangulable by Whitehead's Theorem, by orientability the simplices of the triangulation can be coherently oriented.
The resulting cycle $z$ represents $\left[M\right]$
and clearly satisfies $deg\left(z\right)=1$. Since homologous cycles
have the same degree, a) follows.
b) is part of the definition in \cite[Section 2.1]{ny}.\end{pf}\\
\\
This motivates the following definition.
\begin{df}\label{ifc}Let $M$ be an orientable Riemannian $d$-manifold of nonpositive sectional
curvature, $\widetilde{M}$ its universal cover. We say that $z=\sum_{i=1}^r a_i \sigma_i\in \overline{C}_d^{str}\left(M \right)$ is an ideal fundamental cycle if\\
a) $\partial z=0$, and\\
b) $deg_y\left(z\right)=1$ for all $y\not\in \cup_{i=1}^r\sigma_i\left(\partial\Delta^d\right)$.\\
An ideal fundamental cycle is said to be proper ideal if all $\sigma_i$ are proper.\end{df}

\subsection{Mapping fundamental cycles to ideal fundamental cycles}



\begin{lem}\label{isomorphism4} 
Let $\widetilde{M}$ be a
simply connected Riemannian manifold of nonpositive sectional curvature, $\Gamma$ a discrete group of
isometries of $\widetilde{M}$, $M=\Gamma\backslash\widetilde{M}, \pi:\widetilde{M}\rightarrow M$ the projection, $\tilde{x}_0\in\widetilde{M}, x_0=\pi\left(\tilde{x}_0\right), c_0\in\partial_\infty\widetilde{M}$.

Then there exists a $\Gamma$-invariant chain map $\widetilde{\Psi}:\overline{C}_*^{str}\left(\widetilde{M}\right)
\rightarrow \overline{C}_*^{str}\left(\widetilde{M}\right)$ which is given on 0-simplices by $$\widetilde{\Psi}\left(x\right)=x\mbox{\ if\ }x\in \left(\widetilde{M}-\Gamma\tilde{x}_0\right)\cup\partial_\infty\widetilde{M},$$
$$\widetilde{\Psi}\left(\gamma\tilde{x}_0\right)=\gamma c_0\mbox{\ for\ all\ }\gamma\in\Gamma,$$

such that $\widetilde{\Psi}$ descends to a chain map 
$\Psi:\overline{C}_*^{str,x_0}\left(M\right)\rightarrow \overline{C}_*^{str,c_0}\left(M\right).$ 

If $M$ is a closed, orientable, connected manifold then $\Psi\left(z\right)$ is an ideal fundamental cycle
whenever
$z\in \overline{C}_*^{str,x_0}\left(M\right)$ is an ideal fundamental cycle.
\end{lem}

\begin{pf}
We have defined $\widetilde{\Psi}$ $\Gamma$-equivariantly on the 0-skeleton. Since straight simplices in $\widetilde{M}\cup\partial_\infty\widetilde{M}$ are determined by their vertices, this extends uniquely to a $\Gamma$-equivariant chain map $\widetilde{\Psi}$. Since $\widetilde{\Psi}$ is $\Gamma$-equivariant, it descends to a chain map on $\overline{C}_*^{str}\left(M\right)$. By construction, its restriction to $\overline{C}_*^{str,x_0}\left(M\right)$ has image in $\overline{C}_*^{str,c_0}\left(M\right)$.

To prove the last claim, since $deg$ depends only on the homology class, it suffices to prove $deg\left(\Psi\left(z\right)\right)=1$
for some
ideal fundamental cycle $z$. Choose some fundamental cycle with all vertices in some point $x\not=x_0$, then $\Psi\left(z\right)=z$, which implies $deg\left(\Psi\left(z\right)\right)=1$.

\end{pf}

We remark that $\Psi:\overline{C}_*^{str,x_0}\left(M\right)\rightarrow \overline{C}_*^{str,c_0}\left(M\right)$ is a chain isomorphism if and only if $\Gamma$ acts freely on $\Gamma c_0$, i.e., if and only if $c_0\in\partial_\infty\widetilde{M}$ is not a fixed point of any $\gamma\in\Gamma$.


\subsection{Cusped manifolds}

{\bf Assumption A:} {\em Let $\overline{M}$ be a manifold with boundary components $\partial_1 \overline{M},\ldots, \partial_s \overline{M}$, $M=\overline{M}-\partial \overline{M}$.

Assume that $M$ is a
Riemannian manifold of nonpositive sectional curvature of finite volume and that the universal covering $\widetilde{M}$ is a visibility manifold (\cite[Definition 9.28]{bh}).
Let $\pi:\widetilde{M}\rightarrow M$ be the canonical projection, and let $\widetilde{x_0}\in \widetilde{M}, x_0=\pi\left(\tilde{x}_0\right)$, and $c_0\in\partial_\infty\widetilde{M}$ be fixed.

Let $x_i\in\partial_i\overline{M}$ for $i=1,\ldots,s$ and fix (using a path from $x_0$ to $x_i$) an identification of $\Gamma_i:=\pi_1\left(\partial_i\overline{M},x_i\right)$ with a subgroup of $\Gamma:=\pi_1\left(\overline{M},x_0\right)=\pi_1\left(M,x_0\right)$ (which acts on $\widetilde{M}$).
}\\

Remark: A rank one symmetric space of noncompact type is a visibility manifold. 
If Assumption A holds, then it follows from \cite{ebe} that there are $c_1,\ldots,c_s\in\partial_\infty \widetilde{M}$ 
such that we have a continuous projection
$$p:
\widetilde{M}\bigcup\cup_{i=1}^s
\Gamma c_i\rightarrow DCone\left(\cup_{i=1}^s\partial_i \overline{M}\rightarrow \overline{M}\right) $$
with $x\in\Gamma c_i\Longleftrightarrow
p\left(x\right)$ is the cone point of $Cone\left(\partial_i \overline{M}\right)$, for $i=1,\ldots,s$ (see \cite[Section 4.4]{ku}).

\begin{df}\label{cone} 
If Assumption A holds, then 
a simplex $\sigma\in C_*\left(DCone\left(\cup_{i=1}^s\partial_i \om\rightarrow \om\right)\right)$ is said to be straight if some (hence any) lift $\tilde{\sigma}\in C_*\left(\widetilde{M}\bigcup\cup_{i=1}^s\Gamma c_i\right)\subset C_*\left(\widetilde{M}\cup\partial_\infty\widetilde{M}\right)$ with $p\left(\tilde{\sigma}\right)=\sigma$ is a straight simplex.

If a vertex of $\tilde{\sigma}$ is in $\gamma c_i$ for some $\gamma\in \Gamma, 1\le i\le s$, then we call the corresponding vertex of $\sigma=p\left(\tilde{\sigma}\right)$ an ideal vertex. All other vertices of $\sigma$
are called interior vertices.
Let
$$\widehat{C}_*^{str,x_0}\left(M\right)\subset C_*\left(DCone\left(\cup_{i=1}^s\partial_i \om\rightarrow \om\right)\right)   $$
be the subcomplex freely generated by the straight simplices for which\\
- either all vertices are in $x_0$,\\
- or the last vertex is in $\Gamma c_i$ for some $1\le i\le s$, 
all other vertices are in $x_0$, and the homotopy
classes of all edges
between interior vertices belong to the image of $\Gamma_i$ in $\Gamma$.

$$\widehat{C}_*^{str,c_0}\left(M\right)\subset  C_*\left(DCone\left(\cup_{i=1}^s\partial_i \om\rightarrow \om\right)\right)$$ 
is the subcomplex freely generated by simplices $\tau=\pi\left(\tilde{\tau}\right)$ such that\\
- either all vertices of $\tilde{\tau}$ are in $\Gamma c_0$,\\
- or the last vertex of $\tilde{\tau}$ is in $\Gamma c_i$ for some $1\le i\le s$
and all other vertices of $\tilde{\tau}$ are in $\Gamma c_0$.

\end{df}



\begin{lem}\label{isomorphism3} (\cite{ku}, Lemma 8a)
If Assumption A holds, then
 $$\widehat{C}_*^{str,x_0}\left(M\right)\cong C_*^{simp}\left(DCone\left(\cup_{i=1}^s B\Gamma_i\rightarrow B\Gamma\right)\right).$$
\end{lem}

{\bf Homology invariance of {\em deg} for cusped manifolds.}  If $\om$ is a compact, orientable, connected $d$-manifold with boundary and $z=\sum_{i=1}^ra_i\sigma_i\in
\overline{C}_d^{str}\left(DCone\left(\cup_{i=1}^s\partial_i\overline{M}\rightarrow\overline{M} \right)\right)$ satisfies $\partial z=0$, then for all $y\not\in
\cup_{i=1}^r \sigma_i\left(\partial \Delta^d\right)$ the
isomorphism $$H_d\left(DCone\left(\cup_{i=1}^s\partial_i\overline{M}\rightarrow\overline{M}\right),DCone\left(\cup_{i=1}^s\partial_i\overline{M}\rightarrow\overline{M}\right)-\left\{y\right\}\right)
\rightarrow {\bf Z}$$ sends the relative homology class of $z$
to 
$deg_y\left(z\right)$.
In particular 
$deg_y\left(z\right)$ depends
only on the homology class $\left[z\right]\in H_d\left(DCone\left(\cup_{i=1}^s\partial_i\overline{M}\rightarrow\overline{M}\right)\right)=H_d\left(\om,\partial\om\right)$.

\begin{lem}\label{isomorphism5}
If Assumption A holds, then there is a chain map
$$\hat{\Psi}:\hat{C}_*^{str,x_0}\left(M\right)\rightarrow \hat{C}_*^{str,c_0}\left(M\right),$$
such that the restriction of $\hat{\Psi}$ to $C_*^{str,x_0}$ is the chain map $\Psi$ defined by \hyperref[isomorphism4]{Lemma \ref*{isomorphism4}}.

If $z\in \hat{C}_*^{str,x_0}\left(M\right)$ is an ideal fundamental cycle,
then $\hat{\Psi}\left(z\right)$ is an ideal fundamental cycle.
\end{lem}
\begin{pf}
If $\sigma\in \hat{C}_*^{str,x_0}\left(M\right)$ has all 
vertices in $x_0$, that is if $\sigma\in C_*^{str,x_0}\left(M\right)$, 
then we let $\hat{\Psi}\left(\sigma\right):=\Psi\left(\sigma\right)$, where $\Psi$ is
defined by \hyperref[isomorphism4]{Lemma \ref*{isomorphism4}}. 

In the other case $\sigma$ lifts to a $q$-simplex $\tilde{\sigma}:\Delta^q\rightarrow\widetilde{M}\cup\partial_\infty\widetilde{M}$ whose 
last vertex $v_q$ is $\gamma c_i$ for some $\gamma\in \Gamma, 1\le i\le s$ and we 
have $\partial_q\sigma\in C_{n-1}^{str,x_0}\left(M\right)$. Then
we define $\hat{\Psi}\left(\sigma\right):=p\left(\tilde{\tau}\right)\in C_n^{str,c_0}\left(M\right),$ 
where $\tilde{\tau}$ is the unique straight $q$-simplex with last vertex $\gamma c_0$ and with $\partial_q\tilde{\tau}=\widetilde{\Psi}\left(\partial_q\tilde{\sigma}\right)$.

$\Psi$ is a chain map, hence $\hat{\Psi}\left(\partial\sigma\right)=\partial\hat{\Psi}\left(\sigma\right)$ whenever all vertices of $\sigma$ are interior vertices. In the other case, if $\sigma$ lifts to $\tilde{\sigma}$ with $v_q=\gamma c_i$, then we have $\hat{\Psi}\left(\partial_q\sigma\right)=\Psi\left(\partial_q\sigma\right)=\partial_q
\Psi\left(\sigma\right)=\partial_q\hat{\Psi}\left(\sigma\right)$ and, for $0\le j\le q-1$, $\hat{\Psi}\left(\partial_j\sigma\right)=p\left(\tilde{\kappa}\right)$, where $\tilde{\kappa}$ is the straight $q-1$-simplex with last vertex $\gamma c_0$ and with $\partial_{q-1}\tilde{\kappa}=\widetilde{\Psi}\left(\partial_j\partial_q\tilde{\sigma}\right)=\partial_j\widetilde{\Psi}\left(\partial_q\tilde{\sigma}\right)$, hence $p\left(\tilde{\kappa}\right)=\partial_j\hat{\Psi}\left(\sigma\right)$. Thus $\hat{\Psi}$ is a chain map.


Again, since $deg$ depends only on the homology class, it suffices to prove $deg\left(\hat{\Psi}\left(z\right)\right)=1$
for some
ideal fundamental cycle $z$. Choose some ideal fundamental cycle with all {\em interior} vertices in some point $x\not=x_0$, then $\hat{\Psi}
\left(z\right)=z$, which implies $deg\left(\hat{\Psi}\left(z\right)\right)=1$.

\end{pf}

As a side remark, unrelated to the rest of the paper, 
we observe that a construction analogous to the proof of \hyperref[isomorphism5]{Lemma \ref*{isomorphism5}} 
proves the possibility of an alternative definition of the simplicial volume for hyperbolic manifolds. The simplicial volume of a compact, orientable, connected manifold $\om$, as defined by Gromov, is $\parallel \om,\partial \om\parallel =inf\left\{\sum_{i=1}^r \mid a_i\mid: \sum_{i=1}^ra_i\sigma_i\mbox{\ represents\ }\left[\om,\partial \om\right]\right\}$. 
If we define the ideal simplicial volume as $\parallel M\parallel_{ideal}=inf\left\{\sum_{i=1}^r \mid a_i\mid: \sum_{i=1}^r a_i\tau_i\mbox{\ is\ an\ ideal\ fundamental\ cycle}\right\}$, then we can apply a similar construction to map any relative fundamental cycle to an ideal fundamental cycle, thus proving $\parallel M\parallel_{ideal}
\le \parallel \om,\partial \om\parallel$ whenever $M=\om - \partial\om$ is nonpositively curved. Moreover, if $M$ is
hyperbolic, then each ideal simplex has volume at most $V_n$, thus the ideal simplicial volume can not be bigger than $\frac{1}{V_n}Vol\left(M\right)$. Since $\parallel \om,\partial \om\parallel=\frac{1}{V_n}Vol\left(M\right)$ by the Gromov-Thurston Theorem, this implies 
$$\parallel \om,\partial \om\parallel=\parallel M\parallel_{ideal}$$ for hyperbolic manifolds of finite volume.
\subsection{Well-definedness of $\beta\left(M\right)$}
If $\overline{M}$ is a compact, orientable, connected manifold such that its interior admits a metric of negative sectional curvature, then there exists an ideal fundamental cycle. Indeed, one can take some triangulation $\sum_{i=1}^r\tau_i$ of $\left(\overline{M},\partial\overline{M}\right)$, let $\sum_{j=1}^p\kappa_j:=DCone\left(\partial\left(\sum_{i=1}^r \tau_i\right)\right)$ be the  
cone over the induced triangulations of the path-components $\partial_1 M,\ldots,\partial_sM$ of $\partial M$, with one cone point for each path-component $\partial_i M$, and then use the 
homeomomorphism $DCone\left(\cup_{l=1}^s\partial_i M\rightarrow M\right)\cong \Gamma\backslash G/K \cup \left\{\Gamma c_1,\ldots,\Gamma c_s\right\}$.

Thus we can define $\beta\left(\om\right)$ by \hyperref[nyinvariante]{Definition \ref*{nyinvariante}}.
We will now prove that the definition of $\beta\left(\om\right)$ does not depend on the chosen proper ideal fundamental cycle.

\begin{lem}\label{welldefined}Let $\overline{M}$ be a compact, orientable, connected manifold such that its interior admits a metric of negative sectional curvature and finite volume. Let $G$ be the isometry group of $\widetilde{M}$.
Then $\beta\left(\om\right)$ in \hyperref[nyinvariante]{Definition \ref*{nyinvariante}} is well-defined:
if $\sum_{i=1}^ra_i\tau_i$ and $\sum_{j=1}^sb_j\kappa_j$ are proper ideal fundamental cycles,
then $\sum_{i=1}^ra_i \left[cr\left(\tau_i\right)\right]=\sum_{j=1}^sb_j\left[cr\left(\kappa_j\right)\right]
\in H_*\left(C_*\left(\partial_\infty\widetilde{M}\right)_G\right)$.\end{lem}
\begin{pf} By assumption and \hyperref[ifc]{Definition \ref*{ifc}}, we have $\sum_{i=1}^r a_i deg_x\left(\tau_i\right)=\sum_{j=1}^s b_j deg_x\left(\kappa_j\right)=1$ for all $x\not\in \cup_{i=1}^r\tau_i\left(\partial\Delta^n\right)\cup\cup_{j=1}^s \kappa_j\left(\partial\Delta^n\right)$.

If we define algebraic volume $algvol\left(\sigma\right)$ of straight simplices as in \cite{bp}, p.107, by $algvol\left(\sigma\right)=
\pm vol\left(\sigma\right)$ with the sign according to the orientation of $\sigma$, then\footnotemark\footnotetext[4]{using that $\cup_{i=1}^r\tau_i\left(\partial\Delta^n\right)\cup\cup_{j=1}^s \kappa_j\left(\partial\Delta^n\right)$ is a null set, thus can be neglected for integration } (as in \cite{bp}, p.109):
$$\sum_{i=1}^r a_i algvol\left(\tau_i\right)= 
 \sum_{i=1}^r a_i\int_M \sum_{\tau_i\left(x\right)=y}sign\left(D_x\tau_i\right)
dvol\left(y\right)$$
$$
=\int_M\sum_{i=1}^r a_i deg_y\left(\tau_i\right)dvol\left(y\right)=
\int_M1dvol\left(y\right)=vol\left(M\right),$$
in particular, $\sum_{i=1}^ra_i algvol\left(\tau_i\right)=\sum_{j=1}^s b_jalgvol\left(\kappa_j\right)$. 
$$Algvol: C_*\left(\partial_\infty\widetilde{M}\right)_G\rightarrow{\bf R}$$
$$\left(c_0,\ldots,c_n\right)\rightarrow algvol\left(str\left(c_0,\ldots,c_n\right)\right),$$
where $str\left(c_0,\ldots,c_n\right)$ is the unique proper ideal straight simplex with these vertices, is by Stokes Theorem a (nontrivial) chain map and satisfies of course $Algvol\circ cr=algvol$.

By \hyperref[discrete]{Corollary \ref*{discrete}} we have $H_*\left(C_*\left(\partial_\infty\widetilde{M}\right)_\Gamma\right)\simeq{\bf Z}$, therefore 
$Algvol$ must be injective on $H_*\left(C_*\left(\partial_\infty\widetilde{M}\right)_\Gamma\right)$. 

In particular, $$\sum_{i=1}^ra_i algvol\left(cr\left(\tau_i\right)\right)=\sum_{j=1}^s b_jalgvol\left(cr\left(\kappa_j\right)\right)$$ implies $$\sum_{i=1}^ra_i \left[cr\left(\tau_i\right)\right]=\sum_{j=1}^sb_j\left[cr\left(\kappa_j\right)\right]
\in H_*\left(C_*\left(\partial_\infty\widetilde{M}\right)_\Gamma\right)\simeq{\bf Z},$$ 
hence $$\sum_{i=1}^ra_i \left[cr\left(\tau_i\right)\right]=\sum_{j=1}^sb_j\left[cr\left(\kappa_j\right)\right]
\in H_*\left(C_*\left(\partial_\infty\widetilde{M}\right)_G\right).$$

\end{pf}

\section{Evaluation}

Let $\widetilde{M}=G/K$ be a symmetric space of noncompact type. Fix some $c_0\in\partial_\infty
\widetilde{M}$.
We define the evaluation map $$ev_{G,c_0}:C_*^{simp}\left(BG\right)\rightarrow C_*\left(\partial_\infty \widetilde{M}\right)\ot $$ on generators by $$ev_G\left(g_1,\ldots,g_n\right)=\left(c_0,g_1c_0,g_1g_2c_0,\ldots,g_1g_2\ldots g_nc_0\right)\otimes 1.$$
It is straightforward to check that $ev_{G,c_0}$ extends linearly to a chain map, thus it induces a homomorphism
$$ev_{G,c_0*}:H_*\left(G\right)\rightarrow {\mathcal{P}}_*\left(G/K\right).$$
If $\Gamma_1,\ldots,\Gamma_s$ is a set of subgroups of $G$ such that for $1\le i\le s$ there is some $c_i\in\partial_\infty\widetilde{M}$ with $\Gamma_i\subset Stab\left(c_i\right)$, then we define
$$ev_{G,\Gamma_1,\ldots,\Gamma_s,c_0,c_1,\ldots,c_s}:C_*^{simp}\left(DCone\left(\cup_{i=1}^sB\Gamma_i\rightarrow BG\right)\right)\rightarrow C_*\left(\partial_\infty \widetilde{M}\right)\ot$$ by 
$$ev_{G,\Gamma_1,\ldots,\Gamma_s,c_0,c_1,\ldots,c_s}\left(g_1,\ldots,g_n\right)=ev_G\left(g_1,\ldots,g_n\right)$$ if $\left(g_1,\ldots,g_n\right)\in C_*^{simp}\left(BG\right)$ and by 
$$ev_{G,\Gamma_1,\ldots,\Gamma_s,c_0,c_1,\ldots,c_s}\left(Cone\left(\gamma_1,\ldots,\gamma_{n-1}\right)\right)=\left(c_0,\gamma_1c_0,\gamma_1\gamma_2c_0,
\ldots,\gamma_1\gamma_2\ldots\gamma_{n-1}c_0,c_i\right)\otimes 1$$ if $Cone\left(\gamma_1,\ldots,\gamma_{n-1}\right)\in C_*^{simp}\left(Cone\left(B\Gamma_i\right)\right)$ for some $1\le i\le s$.

The assumption $\Gamma_i\subset Stab\left(c_i\right)$ implies that $ev_{G,\Gamma_1,\ldots,\Gamma_s,c_0,c_1,\ldots,c_s}$ is a chain map, as can be shown by a 
straightforward calculation. Indeed, for $j\not=0$ it 
is immediate from the definition that $\partial_j ev\left(Cone\left(\gamma_1,\ldots,\gamma_{n-1}\right)\right)=ev\left(\partial_j
(Cone\left(\gamma_1,
\ldots,\gamma_{n-1}\right)\right)$, and for $j=0$ we have (omitting the indices of $ev$)
$$\partial_0 ev\left(Cone\left(\gamma_1,\ldots,\gamma_{n-1}\right)\right)=\left(\gamma_1c_0,\gamma_1\gamma_2c_0,\ldots,\gamma_1\ldots\gamma_{n-1}c_0,c_i\right)\otimes 1$$
$$=\left(c_0,\gamma_2c_0,\ldots,\gamma_2\ldots\gamma_{n-1}
c_0,\gamma_1^{-1}c_i\right)\otimes 1=\left(c_0,\gamma_2c_0,\ldots,\gamma_2\ldots\gamma_{n-1}
c_0,c_i\right)\otimes 1 $$
$$=ev\left(Cone\left(\gamma_2,\ldots,\gamma_{n-1}\right)\right)=ev\left(\partial_0 Cone\left(\gamma_1,\ldots,\gamma_{n-1}\right)\right),$$
where the second equality uses the definition of the tensor product (and that $G$ acts trivially on ${\bf Z}$), and the third equality is true because of $\gamma_1\in \Gamma_i\subset Stab\left(c_i\right)$.


\begin{thm}\label{Thm2} Let $\overline{M}$ be a compact, orientable, connected manifold with boundary such that $M:=\overline{M}
-\partial \overline{M}=\Gamma\backslash G/K$ is a finite-volume, locally rank one symmetric space of noncompact type.
Let $\rho:G\rightarrow GL\left(N,{\Bbb C}\right)$ be a representation. For $1\le i\le s$, let $\Gamma_i$ be the subgroup of $\Gamma\simeq\pi_1\left(M,x_0\right)$ corresponding to $\pi_1\left(\partial_iM,x_i\right)$, let $\Gamma_i^\prime:=\rho\left(\Gamma_i\right)$ and fix some $c_0\in\partial_\infty\widetilde{M}$ and $c_i\in Stab\left(\Gamma_i\right)\subset\partial_\infty\widetilde{M}$ for $1\le i\le s$.

Then $$ev_{SL\left(N,{\bf C}\right),\Gamma_1^\prime,\ldots,\Gamma_s^\prime,\rho_\infty\left(c_0\right),\rho_\infty\left(c_1\right),\ldots,\rho_\infty\left(c_s\right)*}\left(\overline{\gamma}\left(\om\right)\right)=\beta_\rho\left(\om\right).$$\end{thm}

\begin{pf}
Since $M$ is a finite-volume locally rank one symmetric space, $\Gamma_i$ must be unipotent. 
Then $\Gamma_i^\prime$ is unipotent and $\gamma\left(M\right)$ is defined by \hyperref[preimage]{Proposition \ref*{preimage}}. Moreover $c_i\in Stab\left(\Gamma_i\right)\subset\partial_\infty\widetilde{M}$ exists by the Remark after Assumption A in Section 3.3. 

Since $G$ is semisimple, $\rho$ is actually a representation $\rho:G\rightarrow SL\left(N,{\bf C}\right)$.

By definition of $ev_{SL\left(N,{\bf C}\right)\Gamma_1^\prime,\ldots,\Gamma_s^\prime}$ the upper square of the diagram commutes, and $\rho$-equivariance of $\rho_\infty$ implies easily that the second square (whose vertical arrows are induced by $B\rho$ resp.\ $\rho_\infty$) commutes.

(For the sake of readability of the diagram we omit the indices $G,\Gamma_1,\ldots,\Gamma_s,c_0,c_1,\ldots,c_s$ resp.\ $SL\left(N,{\bf C}\right),\Gamma_1^\prime,\ldots,\Gamma_s^\prime,\rho_\infty\left(c_0\right),\rho_\infty\left(c_1\right),\ldots,\rho_\infty\left(c_s\right)$ from the notation.)

$$\begin{xy}
\xymatrix{ C_*\left(BSL\left(N,{\bf C}\right)\right)\ar[d]^{j_*}\ar[r]^{ev}&C_*\left(\partial_\infty\left(SL\left(N,{\bf C}\right)/SU\left(N\right)\right)\right)_{SL\left(N,{\bf C}\right)}={\mathcal{P}}_*^N\left({\bf C}\right)\ar[d]^=\\
C_*\left(DCone\left(\cup_{i=1}^s B\Gamma_i^\prime\rightarrow BSL\left(N,{\bf C}\right)\right)\right)\ar[r]^{ev}& C_*\left(\partial_\infty\left(
SL\left(N,{\bf C}\right)/SU\left(N\right)\right)\right)_{SL\left(N,{\bf C}\right)}={\mathcal{P}}_*^N\left({\bf C}\right)\\
C_*\left(DCone\left(\cup_{i=1}^s B\Gamma_i\rightarrow BG\right)\right)\ar[u]^{B\rho}\ar[r]^{ev}&C_*\left(\partial_\infty\widetilde{M}\right)_G={\mathcal{P}}_*\left(\widetilde{M}\right)\ar[u]^{\rho_*}\\
C_*\left(DCone\left(\cup_{i=1}^s B\Gamma_i\rightarrow B\Gamma\right)\right)\ar[r]^{ev}\ar[u]^{i_1}&
C_*\left(\partial_\infty\widetilde{M}\right)_\Gamma\ar[u]^{i_2}\\
\hat{C}_*^{str,x_0}\left(M\right)\ar[r]^{\hat{\Psi}}\ar[u]^{\hat{\Phi}}& \hat{C}_*^{str,c_0}\left(M\right)\ar[u]^{cr}\\
C_*\left(M\cup\left\{\Gamma c_1,\ldots,\Gamma c_s\right\}\right)\ar[u]^{str}&\\
Z_*\left(\overline{M},\partial \overline{M}\right)
\rightarrow C_*\left(DCone\left(\cup_{i=1}^s \partial_i\overline{M}\rightarrow \overline{M}\right)\right)\ar[u]^{\simeq}&.}
 \end{xy}$$
\\
$i_1$ is the obvious inclusion and $i_2=id\otimes id$. Thus the third square commutes. 

$\hat{\Phi}$ is the isomorphism given by \hyperref[isomorphism3]{Lemma \ref*{isomorphism3}}. If $\sigma\in \hat{C}_*^{str,x_0}\left(M\right)$ has all vertices on $x_0$,
then it lifts to a simplex $\tilde{\sigma}$ in $\widetilde{M}$ with vertices $\tilde{x}_0,\gamma_1\tilde{x}_0,
\gamma_1\gamma_2\tilde{x}_0,\ldots,\gamma_1\ldots\gamma_n\tilde{x}_0$ and (by the proof in \cite{ku}) $\hat{\Phi}\left(\tilde{\sigma}\right)=\left(\gamma_1,\gamma_2,\ldots,\gamma_n\right)$.
By
\hyperref[isomorphism4]{Lemma \ref*{isomorphism4}}, $\widetilde{\Psi}\left(\tilde{\sigma}\right)$ is a simplex with vertices $c_0,\gamma_1 c_0,\ldots,\gamma_1\ldots\gamma_nc_0$, thus 
$$cr\left(\hat{\Psi}\left(\sigma\right)\right)=
\left(c_0,\gamma_1 c_0,\ldots,\gamma_1\ldots\gamma_nc_0\right)\otimes_{{\bf Z}\Gamma}1=ev_{G,\Gamma_1,\ldots,\Gamma_s,c_0,c_1,\ldots,c_s}
\left(\gamma_1,\ldots,\gamma_n\right)=ev_{G,\Gamma_1,\ldots,\Gamma_s,c_0,c_1,\ldots,c_s}\left(\hat{\Phi}\left(\sigma\right)\right).$$ 
If $\sigma\in \hat{C}_*^{str,x_0}\left(M\right)$ has its n-th vertex $v_n$ in $\Gamma c_i$,
then it lifts to a simplex $\tilde{\sigma}$ in $\widetilde{M}\cup\partial_\infty\widetilde{M}$ with 
vertices $\tilde{x}_0,\gamma_1\tilde{x}_0,\gamma_1\gamma_2\tilde{x}_0,\ldots,\gamma_1\ldots\gamma_{n-1}\tilde{x}_0,v_n$ and $\hat{\Phi}\left(\sigma\right)$ is the cone over $\left(\gamma_1,\ldots,\gamma_{n-1}\right)$. By
\hyperref[isomorphism4]{Lemma \ref*{isomorphism4}}, $\widetilde{\Psi}\left(\tilde{\sigma}\right)$ is a simplex 
with vertices $c_0,\gamma_1 c_0,\ldots,\gamma_1\ldots\gamma_{n-1}c_0,c_i$, thus 
$$cr\left(\hat{\Psi}\left(\sigma\right)\right)=
\left(c_0,\gamma_1 c_0,\ldots,\gamma_1\ldots\gamma_{n-1}c_0,c_i\right)\otimes_{{\bf Z}\Gamma}1=ev_{G,\Gamma_1,\ldots,\Gamma_s,c_0}\left(Cone\left(\gamma_1,\ldots,\gamma_{n-1}\right)\right)=ev_{G,\Gamma_1,\ldots,\Gamma_s,c_0}\left(\hat{\Phi}\left(\sigma\right)\right).$$
This proves that the 4th square commutes.

$Z_*\left(\overline{M},\partial \overline{M}\right)
\subset C_*\left(M,\partial M\right)$ denotes the group of relative 
cycles. If $z=\sum_{k=1}^r a_k\tau_k\in Z_*\left(\overline{M},\partial 
\overline{M}\right)$, that is if $\partial z=\sum_{j=1}^q b_j\kappa_j\in C_*\left(\partial\overline{M}\right)$, then each $\kappa_j$ has image in some component $\partial_i\om$ of $\partial\overline{M}$ and
$cone\left(z\right)\in  C_*\left(DCone\left(\cup_{i=1}^s \partial_i\overline{M}\rightarrow \overline{M}\right)\right)$ is 
defined as $cone\left(z\right)=\sum_{k=1}^r a_k\tau_k+\sum_{j=1}^q b_j Cone\left(\kappa_j\right)$, where $Cone\left(\kappa_j\right)$ is contained in $Cone\left(\partial_i\om\right)$. If two relative cycles $z_1,z_2$ are relatively homologous, then $cone \left(z_1\right)$
and $cone\left(z_2\right)$ represent the same homology class.

Since $M=\Gamma\backslash G/K$ is a finite-volume locally rank-one symmetric space, there is a homeomorphism $\Gamma\backslash G/K\cup\left\{\Gamma c_1,\ldots,\Gamma c_s\right\}\simeq DCone\left(\cup_{i=1}^s \partial_i\overline{M}\rightarrow \overline{M}\right)$, inducing the isomorphism $C_*\left(DCone\left(\cup_{i=1}^s \partial_i\overline{M}\rightarrow \overline{M}\right)\right)\simeq C_*\left(M\cup\left\{\Gamma c_1,\ldots,\Gamma c_s\right\}\right)$. 

Finally, $str$ is the straightening which homotopes each cycle to a straight cycle with all interior vertices in $x_0$ and all ideal 
vertices remaining fixed during the homotopy. The same argument as in \cite[C.4.3]{bp} (see the proof of Theorem 4 in \cite{ku}) shows that $str$ exists and preserves homology classes.

Now, if $z$ represents $\left[\om,\partial \om\right]$, then $w:=str\left(cone\left(z\right)\right)$ is an ideal fundamental cycle and, by \hyperref[isomorphism5]{Lemma \ref*{isomorphism5}},
$\hat{\Psi}\left(w \right)$ is an ideal fundamental cycle, thus $i_1\left(cr\left(\hat{\Psi}\left(w\right)\right)\right)$
represents $\beta\left(\om\right)$ by \hyperref[nyinvariant]{Definition \ref*{nyinvariante}}.
 Hence (omitting the index of $ev_{SL\left(N,{\bf C}\right),\Gamma_1^\prime,\ldots,\Gamma_s^\prime,\rho_\infty\left(c_0\right)}$ for the sake of readability) $$\beta_\rho\left(\om\right)=\rho_*\beta\left(\om\right)=
\rho_*i_{2*}\left[cr\hat{\Psi}\left(w\right)\right]=
ev_* \left(B\rho\right)_* i_{1*}\left[ \hat{\Phi}\left(w\right)\right]=$$
$$
ev_* \left(B\rho \right)_*i_{1*}\hat{\Phi}_*\left[\om,\partial \om\right]=ev_*\left(B\rho\right)_*i_{1*}\left[\om,\partial\om\right]=
ev_*\overline{\gamma}_\rho
\left(\om\right).$$
Because of $ev_{SL\left(N,{\bf C}\right),\Gamma_1^\prime,\ldots,\Gamma_s^\prime,\rho_\infty\left(c_0\right)}j_*=ev_{SL\left(N,{\bf C}\right)}$ this implies 
$ev\overline{\gamma}_\rho\left(\om\right)=\beta_\rho\left(\om\right)$.
\end{pf}\\
\\

In \cite{su}, A.\ A.\ Suslin constructed a homomorphism $K_3\left({\bf C}\right)\otimes{\bf Q}\rightarrow {\mathcal{B}}\left({\bf C}\right)\otimes{\bf Q}$ which yields an isomorphism 
$K_3^{ind}\left({\bf C}\right)\otimes{\bf Q}\rightarrow {\mathcal{B}}\left({\bf C}\right)\otimes{\bf Q}$.

The following corollary is well known to the experts, at least in the closed case, but seems not to have appeared in written form so far.
\begin{cor}\label{suslin}If $M$ is an orientable, connected, hyperbolic 3-manifold of finite volume, then Suslin's homomorphism 
$$K_3\left({\bf C}\right)\otimes{\bf Q}\rightarrow {\mathcal{B}}\left({\bf C}\right)\otimes{\bf Q}$$ maps the Goncharov invariant $\gamma\left(M\right)$ to the Neumann-Yang invariant $\beta\left(M\right)\otimes 1$.\end{cor}
\begin{pf}
By \cite[Section 5]{ny} we have $\beta\left(M\right)\in{\mathcal{B}}\left({\bf C}\right)\subset {\mathcal{P}}\left({\bf C}\right)$.

\cite{su} constructs a homomorphism $H_3\left(SL\left(3,{\bf C}\right);{\bf Z}\right)\rightarrow {\mathcal{B}}\left({\bf C}\right)$ such that
restriction to the image of
$H_3\left(SL\left(2,{\bf C}\right);{\bf Z}\right)$ gives our evaluation map $ev_{SL\left(2,{\bf C}\right)}$ (cf.\ \cite[Lemma 3.4]{su}). Composition with the inverse of the stability isomorphism
$H_3\left(SL\left(3,{\bf C}\right);{\bf Z}\right)\rightarrow
H_3\left(SL\left({\bf C}\right);{\bf Z}\right)$ yields a homomorphism $H_3\left(SL\left({\bf C}\right);{\bf Z}\right)\rightarrow {\mathcal{B}}\left({\bf C}\right)$ whose restriction
to the image of
$H_3\left(SL\left(2,{\bf C}\right);{\bf Z}\right)$ gives $ev_{SL\left(2,{\bf C}\right)}$, in particular by \hyperref[Thm2]{Theorem \ref*{Thm2}} it 
maps $\overline{\gamma}\left(M\right)$ to $\beta\left(M\right)\otimes 1$.

Further composition with the map from \cite[Corollary 5.2]{su} gives the homomorphism $K_3\left({\bf C}\right)\rightarrow
{\mathcal{B}}\left({\bf C}\right)$. 
By \cite[Theorem 5.1]{su} 
it induces an isomorphism $K_3\left({\bf C}\right)/\pi_3\left(BGM\left({\bf C}\right)^+\right)\rightarrow {\mathcal{B}}\left({\bf C}\right)$. Indeed, \cite[Theorem 5.1]{su} gives an isomorphism
$K_3\left(F\right)/\pi_3\left(BGM\left(F\right)^+\right)\rightarrow {\mathcal{B}}\left(F\right)/2c$ for any field $F$, where $c\in{\mathcal{B}}\left(F\right)$ is an element of order 6 by \cite[Lemma 1.4]{su}. In particular $c=0$ since ${\mathcal{B}}\left({\bf C}\right)$ is a ${\bf Q}$-vector space.
It follows then from \cite[Theorem 5.2]{su} that this morphism induces the isomorphism $K_3^{ind}\left({\bf C}\right)\otimes{\bf Q}\rightarrow {\mathcal{B}}\left({\bf C}\right)\otimes{\bf Q}$.

		We note that the inclusion $K_3\left({\bf C}\right)\otimes{\bf Q}\cong
PH_3\left(GL\left({\bf C}\right);{\bf Q}\right)\cong PH_3\left(SL\left({\bf C}\right);{\bf Q}\right)\subset H_3\left(SL\left({\bf C}\right);{\bf Q}\right)$ is actually an equality, hence (cf.\ \cite[Section 2.5]{ku} and the definition of $\gamma\left(M\right)$ in \cite[Theorem 4]{ku}) 
$\overline{\gamma}\left(M\right)=\gamma\left(M\right)$. And the claim follows.
\end{pf}

\noindent
Thilo Kuessner\\
Hensenstra\ss e 160\\
D-48149 M\"unster\\
Germany\\
e-mail: kuessner@math.uni-muenster.de\\


\begin{thebibliography}{29}

\bibitem{bp}R.\ Benedetti, C.\ Petronio, 'Lectures on
Hyperbolic Geometry', Universitext, Springer-Verlag, Berlin (1992).
\bibitem{bh}M.\ Bridson, A.\ Haefliger, 'Metric spaces of non-positive curvature', Grundlehren der Mathematischen Wissenschaften 319, Springer-Verlag, Berlin (1999).

\bibitem{dd}W.\ Dicks, M.\ J.\ Dunwoody, 'Groups acting on graphs'. Camb.\ Studies in Adv.\ Math.\ {\bf 17} (1989).
\bibitem{ebe}P.\ Eberlein, 'Lattices in spaces of nonpositive curvature', {\em Ann.\ Math.\ } 111, pp.435-476 (1980).


\bibitem{gon}A.\ Goncharov, 'Volumes of hyperbolic manifolds
and mixed Tate motives', {\em J.\ AMS} 12, pp.569-618 (1999).
\bibitem{ku}T.\ Kuessner, 'Locally symmetric spaces and K-theory of number fields', http://xxx.uni-augsburg.de/abs/0904.0919

\bibitem{ny}W.\ Neumann, J.\ Yang, 'Bloch invariants of hyperbolic 3-manifolds', {\em Duke Math.\ J.} 96, pp.29-59 (1999).
\bibitem{su}A.\ A.\ Suslin, '$K_3$ of a field, and the Bloch group', {\em Trudy Mat.\ Inst.\ Steklov} 183, pp.180-199 (1990).

\end{thebibliography}
\end{document}